\documentclass[12pt, leqno, a4paper]{article}

\usepackage[top=3cm, bottom=4cm, left=3cm, right=3cm]{geometry}

\usepackage{latexsym,amssymb}
\usepackage{amsmath,amsthm}
\usepackage{color} 
\usepackage{mathtools}

\usepackage{multirow,eepic} 

\newtheorem{Theorem}{\bf Theorem}[section]
\newtheorem{Lemma}{\bf Lemma}[section]
\newtheorem{Proposition}{\bf Proposition}[section]
\newtheorem{Corollary}{\bf Corollary}[section]
\newtheorem{Remark}{\bf Remark}[section]

\newenvironment{prf}[1]
   {{\noindent \bf Proof of {#1}.}}{\hfill \qed}

\allowdisplaybreaks[4]
\numberwithin{equation}{section}

\begin{document}
\title{
Attainability of the best Sobolev constant
in a ball\\
} 
\author{
        Norisuke Ioku\footnote{e-mail address: ioku@ehime-u.ac.jp} \\ \\ 
        {\small Graduate School of Science and Engineering, Ehime University} \\ 
        {\small Matsuyama, Ehime 790-8577, Japan} 
        } 
\date{} 
\pagestyle{myheadings}

\maketitle 
\begin{abstract} 
The best constant of the Sobolev inequality in
the whole space is attained by the Aubin--Talenti function; 
however,
this does not happen 
in bounded
domains because the break in dilation invariance.
In this paper, we investigate a new 
scale invariant form of the Sobolev inequality in a ball 
and show that its best constant is attained by functions of the Aubin--Talenti type.
Generalization to the Caffarelli--Kohn--Nirenberg inequality 
in a ball is also discussed.
\end{abstract} 
\noindent 
{\small 
\vspace{5pt} 
}
 

\section{Introduction and Results} 
\label{section:1} 
%

%
The Sobolev inequality
states that,
if $n\ge 2$ and $1\le p<n$, then
\begin{equation}\label{eq:1.1a}
 S_{n,p}\|u\|_{L^{p^*}}
 \le \|\nabla u\|_{L^p}
\end{equation}
for every $u\in W^{1,p}(\mathbb{R}^n)$,
where
$p^*$ is the Sobolev conjugate number defined by $p^*=np/(n-p)$
and 
$S_{n,p}$ is the best constant for the inequality \eqref{eq:1.1a}  given by
\[
\left\{
\begin{aligned}
 S_{n,p}
&
 =
 \sqrt{\pi}n^{\frac1p}\left(\frac{n-p}{p-1}\right)^{\frac{p-1}{p}}
 \left[
  \frac{
    \Gamma\left(\frac{n}{p}\right)\Gamma\left(n+1-\frac{n}{p}\right)
   }
   {
   \Gamma(n)\Gamma\left(1+\frac{n}{2}\right)
   }
 \right]^{\frac1n}\ \ 
 &&
 \text{for}\ 1<p<n,
 \\
 S_{n,1}
 &
 =\sqrt{\pi}\frac{n}{\left[\Gamma\left(1+\frac{n}{2}\right)\right]^{\frac{1}{n}}}
 &&
 \text{for}\ p=1,
\end{aligned}
\right.
\]
where $\Gamma(\cdot)$ is the gamma function.
The inequality \eqref{eq:1.1a} with the best constant
was obtained by 
Federer-Fleming~\cite{FF} and Maz$'$ya~\cite{Ma1}
for $p = 1$ and by Aubin~\cite{Au} and Talenti~\cite{Ta1} for $1 < p < n$.
Furthermore, if $1 < p < n$, the best constant is attained by the two parameter family 
\begin{equation}\label{AubinTalenti}
U(x)=(a+b|x|^{\frac{p}{p-1}})^{1-\frac{n}{p}},\quad a,b>0
\end{equation}
and its translation. 
If we replace $\mathbb{R}^n$ in $\Omega \subset \mathbb{R}^n$, the Sobolev inequality still holds for $u\in W^{1,p}_0(\Omega)$ with the same best constant;
however, it is not attained in $W^{1,p}_0(\Omega)$
because the dilation invariance 
under
\begin{equation}\label{eq:1.7a}
u_{\mu}(x)=\mu^{\frac{n-p}{p}}u(\mu x),\qquad \mu>0,
\end{equation}
breaks
(see \cite[p.44]{Ke} or \cite[Proposition~1.43]{Wi}). 


On the other hand, {in the critical case of $p=n$}, 
{the Sobolev inequality of the form}  
\eqref{eq:1.1a} 
is no longer true because $S_{n,p}\to 0$ as $p\to n$, and 
$W^{1,n}_0(\Omega)$ is not embedded into $L^{\infty}(\Omega)$.
Alvino~\cite{A} considered the critical case and obtained
\begin{equation}\label{eq:1.6aa}
\sqrt{\pi}\frac{n^{\frac1n}}{\Gamma\left(1+\frac{n}{2}\right)^{\frac1n}}
 \sup_{x\in B_R}\frac{|u^{*}(x)|}{\left(\log\left(\frac{R}{|x|}\right)\right)^{\frac{n-1}{n}}}
 \le 
 \|\nabla u\|_{L^n(B_R)}\qquad \text{for\ all}\ u\in W^{1,n}_0(B_R),
\end{equation}
where $B_R$ is the  ball centered at the origin with radius $R>0$
and 
$u^{*}$ denotes the Schwarz symmetrization of $u\in W^{1,n}_0(B_R)$. For the definition of $u^{*}$, see the end of Section~\ref{section:1} of this paper or \cite[Section~1.3]{Ke}.
Alvino's inequality~\eqref{eq:1.6aa} is known to be the critical case of Sobolev embeddings because
this inequality implies the optimal embedding of $W^{1,n}_0(B_R)$ into Orlicz spaces
(see~\cite[Example~1]{C}).
Several equivalent forms of \eqref{eq:1.6aa}, especially the relationship between \eqref{eq:1.6aa} and Moser--Trudinger inequalities, are discussed by Cassani--Sani--Tarsi~\cite{CST}. 


One of the main difference between the Alvino inequality~\eqref{eq:1.6aa} and the Sobolev inequality~\eqref{eq:1.1a}
is the scale invariance structure, that is,
the critical case~\eqref{eq:1.6aa} is not invariant under the dilation $x\mapsto \lambda x$ 
but is invariant 
under 
\begin{equation}\label{eq:1.8a}
 u_{\lambda}(x)=\lambda^{-\frac{n-1}{n}}u\left(\left(\frac{|x|}{R}\right)^{\lambda-1}x\right),\qquad 
 \lambda>0.
\end{equation}
This scaling was firstly found by Adimurthi--do \'O--Tintarev~\cite{AOT},
and they pointed out that Moser functions are invariant under \eqref{eq:1.8a}.
Cassani--Ruf--Tarsi~\cite{CRT} proved that Moser functions are the minimizers of the minimizing problem associated with \eqref{eq:1.6aa} by focusing its invariance under \eqref{eq:1.8a}.
Costa--Tintarev~\cite{CT} applied the scaling \eqref{eq:1.8a} to analyze the concentration profiles of the Trudinger--Moser functional.


\bigskip

It is well known that the inequalities~\eqref{eq:1.1a} and \eqref{eq:1.6aa} with the best constant 
have vast applications and generalizations
in geometry, physics, and functional analysis.
Even though 
\eqref{eq:1.1a} and \eqref{eq:1.6aa} are widely studied, several questions still arise naturally
from the view-point of attainability of the best constant and the scale invariance property.
Does the Sobolev inequality in $W^{1,p}_0(B_R)$ have a scale invariant form?
Can we obtain Alvino's inequality~\eqref{eq:1.6aa} from the Sobolev type inequalities by directly limiting the procedure $p\to n$?
Is there any relationship between the two scalings: the dilation \eqref{eq:1.7a} for $W^{1,p}(\mathbb{R}^n)$ 
and the scaling \eqref{eq:1.8a} for $W^{1,n}_0(B_R)$?

\bigskip

In this paper,
we 
give a positive answer for these questions
by discovering
the scale invariant form of a Sobolev type inequality in $W^{1,p}_0(B_R)$,
which recovers the attainability of the best constant and implies Alvino's inequality by taking $p\to n$.

To state our result, we introduce 
the $q$-logarithmic function and $q$-exponential function as follows:
\begin{equation}\label{eq:2.1a}
\begin{aligned}
\log_q r
&
:=\frac{r^{1-q}-1}{1-q},
\\
\exp_q(r)
&
:=[1+(1-q)r]^{\frac{1}{1-q}},
\end{aligned}
\end{equation}
for $q>0,\ q\neq 1$ and $r>0$.
It is easy to verify that 
\[
\lim_{q\to 1}\log_q r=\log r,\qquad \lim_{q\to 1}\exp_q r=e^r\qquad \text{for\ all}\ r>0.
\]
These modified logarithmic, exponential functions were originally introduced by
Tsallis~\cite{Ts1}
to study nonextensive statistics.
See \cite[Section~3]{Ts} or \cite[Section~2 and Appendix~A]{Suyari} for more details on $q$-logarithmic, exponential functions.

The first result is an improved Sobolev inequality for radially symmetric functions.

%


\begin{Theorem}[Sobolev type inequality]\label{Theorem:1.2}
Let $n\in \mathbb{N},\ n\ge 2$, $1<p<n$,
and
$\frac{1}{p^*}=\frac1p-\frac1n$.
Then for any
radially symmetric
function
$u\in W^{1,p}_0(B_R)$
the following holds:
\begin{equation}\label{eq:1.4}
 S_{n,p}
 \left(\frac{n-p}{p-1}\right)^{\frac1n-1}
 \left(
 \int_{B_{R}}
 \frac{
 |u(x)|^{p^{*}}
 }
 {
 \left[
 \log_{{\frac{n-1}{p-1}}}\frac{R}{|x|}
 \right]^{\frac{p(n-1)}{n-p}}
 }
 dx
 \right)^{\frac{1}{p^{*}}}
\le
\|\nabla u\|_{L^p(B_R)}.
\end{equation}
The left-hand-side constant is optimal and attained by 
\begin{equation}\label{eq:1.9a}
\begin{aligned}
U_R(x)
=
\left[a+b
\left\{
\frac{1}{|x|^{\frac{n-p}{p-1}}}
-
\frac{1}{R^{\frac{n-p}{p-1}}}
\right\}^{-\frac{p}{n-p}}\right]^{1-\frac{n}{p}}
\end{aligned}
\end{equation}
where
$a,b>0$.
Furthermore,
the inequality~\eqref{eq:1.4}
is invariant under the following nonlinear scaling
\begin{equation}\label{eq:1.3}
\left\{
\begin{aligned}
 &
 u_{\lambda}(x)
 := \lambda^{-\frac{p-1}{p}}u(x_{\lambda}),
 \\
 &
 x_{\lambda}
 :=
 \left[
  \lambda |x|^{-\frac{n-p}{p-1}}+(1-\lambda)R^{-\frac{n-p}{p-1}}
 \right]^{-\frac{p-1}{n-p}}\frac{x}{|x|}.
\end{aligned}
\right.
\end{equation}
\end{Theorem}
\begin{Remark}
The inequality~\eqref{eq:1.4} yields the classical Sobolev inequality in a bounded domain $\Omega$ and non-attainability of the best constant by symmetrization techniques as follows. 
Let $u\in W^{1,p}_0(\Omega)$ and $u^{*}$ be its Schwarz symmetrization, where $B_R$ is the ball centered at the origin having the same measure as $\Omega$.
Since
\begin{equation}\label{eq:1.123ccc}
 \frac{n-p}{p-1}
 \log_{{\frac{n-1}{p-1}}}\frac{R}{|x|}
 =
 1-\left(\frac{|x|}{R}\right)^{\frac{n-p}{p-1}}
 < 1
\end{equation}
for all $x\in B_R$,
the inequality 
\eqref{eq:1.4}
to $u^{*}\in W^{1,p}_0(B_R)$ 
together with the P\'olya-Szeg\"o inequality
shows that
\begin{equation}
\label{eq:1.5a}
\begin{aligned}
S_{n,p}
\|u\|_{L^{p^*}(\Omega)}
=
S_{n,p}
\|u^{*}\|_{L^{p^*}(B_R)}
&
<
 S_{n,p}
 \Bigg(
 \int_{B_{R}}\frac{|u^{*}(x)|^{p^{*}}}
 {
 \left[
 \frac{n-p}{p-1}
 \log_{{\frac{n-1}{p-1}}}\frac{R}{|x|}
 \right]^{\frac{p(n-1)}{n-p}}
 }dx
 \Bigg)^{\frac{1}{p^{*}}}
 \\
&
\le
\|\nabla u^{*}\|_{L^{p}(B_R)}
\le
\|\nabla u\|_{L^{p}(\Omega)}
\end{aligned}
\end{equation}
for all $u\in W^{1,p}_0(\Omega)$.
Moreover, the strict inequality~\eqref{eq:1.5a} directly shows that there is no extremal function
for the classical Sobolev inequality in the ball $B_R$.
\end{Remark}

\begin{Remark}
Clearly the equality in \eqref{eq:1.123ccc} shows that
Theorem~\ref{Theorem:1.2}
yields 
the classical Sobolev inequality in $\mathbb{R}^n$
by taking $R\to \infty$.
The scaling~\eqref{eq:1.8a} converges to 
the dilation~\eqref{eq:1.7a}
as $R\to \infty$, that is,
\[
 \lambda^{-\frac{p-1}{p}}u(x_{\lambda})
 \to 
 (\lambda^{-\frac{p-1}{n-p}})^{\frac{n-p}{p}}u\left(\lambda^{-\frac{p-1}{n-p}} x\right),\quad R\to \infty.
\]
Furthermore, the extremal function $U_R(x)$ in \eqref{eq:1.9a} converges pointwise to the Aubin--Talenti function $U(x)$ in \eqref{AubinTalenti} as $R\to \infty$.
\end{Remark}

\begin{Remark}
While the classical Sobolev inequality~\eqref{eq:1.1a}
does not imply
Alvino's inequality, 
%
Theorem~\ref{Theorem:1.2} 
yields Alvino's inequality 
by the direct limiting procedure $p\to n$.
Indeed, 
the explicit value of $S_{n,p}$
gives us that
\[
\begin{aligned}
S_{n,p}\left(\frac{n-p}{p-1}\right)^{-\frac{n-1}{n}}
\to
\frac{\sqrt{\pi}n^{\frac1n}}{\Gamma\left(1+\frac{n}{2}\right)^{\frac1n}}
\end{aligned}
\]
as $p\to n$. 
This together with
$
{\displaystyle \lim_{p\uparrow n}}
\log_{\frac{n-1}{p-1}}\frac{R}{|x|}
=
\log\frac{R}{|x|}
$
yields the desired convergence
\[
S_{n,p}\left(\frac{n-p}{p-1}\right)^{-\frac{n-1}{n}}
\left(
 \int_{B_{R}}\frac{|u^{*}(x)|^{p^{*}}}
 {
 \left[\log_{{\frac{n-1}{p-1}}}\frac{R}{|x|}\right]^{\frac{p(n-1)}{n-p}}
 }dx
 \right)^{\frac{1}{p^{*}}}
\to
\frac{\sqrt{\pi}n^{\frac1n}}{\Gamma\left(1+\frac{n}{2}\right)^{\frac1n}}
  \sup_{x\in B_R}\frac{|u^{*}(x)|}{\left(\log
  \frac{R}{|x|}
  \right)^{\frac{n-1}{n}}}
\]
as $p\to n$.
Furthermore, 
the scaling~\eqref{eq:1.3} coincides with
the scaling~\eqref{eq:1.8a} if $p\to n$,
since
$
{\displaystyle \lim_{p\to n}}x_{\lambda}=\left(\frac{|x|}{R}\right)^{\lambda-1}x.
$

\end{Remark}

Theorem~\ref{Theorem:1.2} does not hold for arbitrary (not necessarily radially symmetric) $W^{1,p}_0(B_R)$ functions. To prove this, 
let us assume that \eqref{eq:1.4} holds for arbitrary $W^{1,p}_0(B_R)$ functions and then derive a contradiction. 
Let $d(x)$ be the distance function to the boundary defined by $d(x):=R-|x|$.
Since there exists $C>0$ such that $d(x)/C\le \log_{\frac{n-1}{p-1}}\frac{R}{|x|}\le Cd(x)$
for all $x\in B_R$,
the inequality \eqref{eq:1.4} yields the following Hardy--Sobolev inequality with the weight function:
\begin{equation}\label{eq:1.10a}
\displaystyle
 \left(
 \int_{B_{R}}|u(x)|^{p^{*}}d(x)^{-\frac{p(n-1)}{n-p}}dx
 \right)^{\frac{1}{p^{*}}}
 \le
 C\|\nabla u\|_{L^p(B_R)}.
\end{equation}
The necessary condition for the exponent of the weight function is known in the study of Hardy--Sobolev inequalities with weights(see \cite[Theorem~19.10 and Remark~19.13 with $\kappa=1,q=0$]{OK}).
Since $-\frac{p(n-1)}{n-p}$ is not admissible, the inequality \eqref{eq:1.10a} cannot hold.
Nevertheless, we obtain the following Sobolev type inequality for arbitrary $W^{1,p}_0(B_R)$ functions
by introducing the differential operator $L_{p}$.
The radial derivative $\nabla_r u(x)$ 
and
the tangential derivative 
$\nabla_{\mathbb{S}^{n-1}}u(x)$ are defined by
\begin{equation}
\label{eq:1.30c}
\nabla_ru(x):=\left(\frac{x}{|x|}\cdot \nabla u(x)\right)\frac{x}{|x|},\qquad 
\nabla_{\mathbb{S}^{n-1}}u(x):=\nabla u(x)-\nabla_r u(x).
\end{equation}
We define $L_p$ as follows:
\begin{equation*}
L_pu:=
\nabla_{\mathbb{S}^{n-1}}u(x)
\slash
\left[
 \frac{n-p}{p-1}
 \log_{{\frac{n-1}{p-1}}}\frac{R}{|x|}
\right]
+
\nabla_r u(x).
\end{equation*}




Clearly the following theorem includes 
Theorem~\ref{Theorem:1.2}, 
since 
$
L_pu=\nabla_r u=\nabla u
$
if $u$ is radially symmetric.

\begin{Theorem}\label{Theorem:1.3}
Let $n\in \mathbb{N},\ n\ge 2$, $1<p<n$.
It follows that
\begin{equation}\label{eq:1.2}
\begin{aligned}
S_{n,p}\left(\frac{n-p}{p-1}\right)^{-\frac{n-1}{n}}
\left(
 \int_{B_{R}}\frac{|u(x)|^{p^{*}}}
 {
 \left[\log_{{\frac{n-1}{p-1}}}\frac{R}{|x|}\right]^{\frac{p(n-1)}{n-p}}
 }dx
 \right)^{\frac{1}{p^{*}}}
\le
\Biggl(
 \int_{{B_{R}}}
|L_p u(x)|^{p}
 dx
\Biggr)^{{\frac{1}{p}}}
\end{aligned}
\end{equation}
for all $u\in C_0^{\infty}(B_R)$.
The constant in the left hand side is optimal and attained by 
$U_R(x)$ defined in \eqref{eq:1.9a}.
Furthermore,
the inequality~\eqref{eq:1.2}
is invariant under the scaling~\eqref{eq:1.3}.
\end{Theorem}

\bigskip

%
One can further generalize \eqref{Theorem:1.3} to inequalities of the Caffarelli--Kohn--Nirenberg type,
which are weighted versions of Sobolev inequalities (See \cite{CKN,HK}).
The Caffarelli--Kohn--Nirenberg inequality is stated as follows.
For any $1< p<n,\ \theta> 1,\ 0\le \frac{1}{\theta}-\frac{1}{\sigma}\le \frac{1}{n}$
it follows that for all $u\in C_0^{\infty}(\mathbb{R}^n)$, we have
\begin{equation}\label{eq:1.5bbb}
S
\left(
 \int_{\mathbb{R}^n}|x|^{\frac{n\sigma}{p^*}}|u(x)|^{\sigma}\frac{dx}{|x|^n}
\right)^{\frac{1}{\sigma}}
\le
\left(
 \int_{\mathbb{R}^n}|x|^{\frac{n\theta}{p}}|\nabla u(x)|^{\theta}\frac{dx}{|x|^n}
\right)^{\frac{1}{\theta}}.
\end{equation}
Furthermore, if $\theta=\sigma$, it follows that for all not necessarily radial functions $u\in C_0^{\infty}(\mathbb{R}^n)$, we have
\begin{equation}\label{eq:1.5ccc}
\tilde{S}
 \int_{\mathbb{R}^n}|x|^{\frac{n\theta}{p^*}}|u(x)|^{\theta}\frac{dx}{|x|^n}
\le
 \int_{\mathbb{R}^n}|x|^{\frac{n\theta}{p}}|\nabla_r u(x)|^{\theta}\frac{dx}{|x|^n}.
\end{equation}
This homogeneous case is an improvement of \eqref{eq:1.5bbb} in the sense that $|\nabla_r u(x)|\le |\nabla u(x)|$ for all $x\in \mathbb{R}^n$.
Here $S$ and 
$\tilde S$ denote the best constant of \eqref{eq:1.5bbb} and \eqref{eq:1.5ccc}, respectively.
The Caffarelli--Kohn--Nirenberg inequality~\eqref{eq:2.4} 
coincides with the Sobolev inequality if $\theta=p$ and $\sigma=p^*$.
The same inequality holds on a bounded domain $\Omega$; however, the dilation invariance $x\mapsto \lambda x$ is broken.
Our main result corresponds with the scale invariance form of \eqref{eq:1.5bbb} in a ball.
%
%
Let 
$1<p<n,\ \theta> 1,\ 0\le \frac{1}{\theta}-\frac{1}{\sigma}\le \frac{1}{n}$,
$\frac{1}{p^{*}}=\frac{1}{p}-\frac{1}{n}$,
and
\begin{equation}\label{eq:2.30c}
q=1+\frac{n-p}{p}\frac{\theta}{\theta-1}.
\end{equation}
For the main result, the differential operator $L_{p,\theta}$ is defined by
\begin{equation}\label{eq:1.15a}
L_{p,\theta}u
:=
\nabla_{\mathbb{S}^{n-1}}u(x)
\slash
\left[
(q-1)
\log_{q}\frac{R}{|x|}
\right]
+
\nabla_r u(x).
\end{equation}
This generalize $L_{p}$ in Theorem~\ref{Theorem:1.3} since $L_{p,p}=L_p$.
%
Our main result is given below as Theorem~\ref{Theorem:1.1}.
\begin{Theorem}[Inequality of the Caffarelli--Kohn--Nirenberg type in $B_R$]\label{Theorem:1.1}
Let $1< p<n,\ \theta> 1,\ 0\le \frac{1}{\theta}-\frac{1}{\sigma}\le \frac{1}{n}$, 
$\frac{1}{p^{*}}=\frac{1}{p}-\frac{1}{n}$, and
$q=1+\frac{n-p}{p}\frac{\theta}{\theta-1}$.
Then, it follows that 
%
%
\begin{equation}\label{eq:3.17a}
\begin{aligned}
S
\left(
 \int_{{B_{R}}}\frac{|x|^{\frac{n\sigma}{p^{*}}}|u(x)|^{{\sigma}}
 }
 {
\left[
(q-1)
\log_{q}\frac{R}{|x|}
\right]^{{1+\frac{\theta-1}{\theta}\sigma}}
 }
 \frac{dx}{|x|^{{n}}}
\right)^{{\frac{1}{\sigma}}}
\le
\left(
 \int_{{B_{R}}}
 |x|^{\frac{n\theta}{p}}
 |L_{p,\theta}u|^{\theta}
 \frac{dx}{|x|^{{n}}}
\right)^{{\frac{1}{\theta}}}
\end{aligned}
\end{equation}
for all $u\in C_0^{\infty}(B_R)$ with the same constant $S$ as given in \eqref{eq:1.5bbb}.
If in addition $\theta=\sigma$, it follows
\begin{equation}\label{eq:3.17ccc}
\begin{aligned}
\tilde{S}
 \int_{{B_{R}}}\frac{|x|^{\frac{n\theta}{p^{*}}}|u(x)|^{{\theta}}
 }
 {
\left[
(q-1)
\log_{q}\frac{R}{|x|}
\right]^{\theta}
 }
 \frac{dx}{|x|^{{n}}}
\le
 \int_{{B_{R}}}
 |x|^{\frac{n\theta}{p}}
 |
\nabla_r u|^{\theta}
 \frac{dx}{|x|^{{n}}}
\end{aligned}
\end{equation}
for all $u\in C_0^{\infty}(B_R)$ with the same constant $\tilde S$ as given in \eqref{eq:1.5ccc}.

Moreover, the constants $S$ and $\tilde S$ 
are optimal, and \eqref{eq:3.17a} and \eqref{eq:3.17ccc} are invariant under 
\begin{equation}\label{eq:2.1a}
\left\{
\begin{aligned}
 &
 u_{\lambda}(x)
 := \lambda^{-\frac{\theta-1}{\theta}}u(x_{\lambda}),
 \\
 &
 x_{\lambda}
 :=
 \left[
  \lambda |x|^{-\frac{n-p}{p}\frac{\theta}{\theta-1}}+(1-\lambda)R^{-\frac{n-p}{p}\frac{\theta}{\theta-1}}
 \right]^{-\frac{p}{n-p}\frac{\theta-1}{\theta}}\frac{x}{|x|}.
\end{aligned}
\right.
\end{equation}

\end{Theorem}

%
The Sobolev type inequalities \eqref{eq:1.2} and \eqref{eq:1.4} can be obtained as a direct consequence from \eqref{eq:3.17a} by taking 
$\theta=p,\ \sigma=p^*$.
%

%
If $\theta=\sigma=p$, the inequality \eqref{eq:1.5ccc} is known as Hardy's inequality and is given as follows:
\[
\frac{n-p}{p}
\int_{\mathbb{R}^n}\frac{|v(x)|^p}{|x|^p}
\le
\int_{\mathbb{R}^n}|\nabla_r v(x)|^pdx,\quad v\in W^{1,p}(\mathbb{R}^n),
\]
where the constant is the best possible and the inequality is invariant under the dilation $x\mapsto \lambda x$.
Theorem~\ref{Theorem:1.1} gives us a corresponding result which has scale invariance in a ball.
\begin{Corollary}
Let $1<p<n$. There holds 
\begin{equation}\label{eq:2.10a}
\left(\frac{p-1}{p}\right)^p
\int_{B_R}\frac{|u(x)|^p}{|x|^p\left[\log_{\frac{n-1}{p-1}}\frac{R}{|x|}\right]^p}dx
\le
\int_{B_R}|\nabla_r u(x)|^p dx,\quad u\in W^{1,p}_0(B_R),
\end{equation}
where the constant is optimal and the inequality is invariant under the scaling \eqref{eq:1.3}.
\end{Corollary}
To the best knowledge of the author, the elemental inequality \eqref{eq:2.10a} is new and natural
because \eqref{eq:2.10a} derives the critical Hardy inequality 
\[
\left(\frac{n-1}{n}\right)^n
\int_{B_R}\frac{|u(x)|^n}{|x|^n\left(\log\frac{R}{|x|}\right)^n}dx 
\le 
\int_{B_R}|\nabla_r u(x)|^ndx,\quad u\in W^{1,n}_0(B_R)
\]
by taking $p\to n$.
For more details on Hardy's inequality, see \cite{
Da,II1,II2,ST,Takahashi} and the references therein.

\bigskip

We now summarize some notations and basic facts used in this paper.
Let $B_R$ be the ball centered at the origin with radius $R>0$.
The surface measure of the unit sphere $\mathbb{S}^{n-1}$ is denoted by $\omega_{n-1}$.
The Schwarz symmetrization of $f$ denoted by $f^{*}$
is defined by 
$$
   f^*(x)
=\inf\left\{\lambda>0:\mu_{f}(\lambda)\le \frac{\omega_{n-1}}{n}|x|^n\right\},
$$
where $\mu_f(\lambda)$ is the distribution function of $f$.
The Schwarz symmetrization $f^{*}$ is radially symmetric and non-increasing.
For $f\in W^{1,p}(\mathbb{R}^n)$, the Schwarz symmetrization $f^*$ belongs to $W^{1,p}(\mathbb{R}^n)$ and 
the P\'olya-Szeg\"o inequality
\[
\|\nabla f^*\|_{L^p(\mathbb{R}^n)}\le \|\nabla f\|_{L^p(\mathbb{R}^n)}
\]
holds. For more details on the Schwarz symmetrization, see \cite{Ke}.
The $q$-logarithmic function and $q$-exponential function are defined 
by
\begin{equation*}
\begin{aligned}
\log_q r
:=\frac{r^{1-q}-1}{1-q},
\qquad 
\exp_q(r)
:=[1+(1-q)r]^{\frac{1}{1-q}},
\end{aligned}
\end{equation*}
for $q>0,\ q\neq 1$ and $r>0$.
It is easy to verify that 
\[
\lim_{q\to 1}\log_q r=\log r,\qquad \lim_{q\to 1}\exp_q r=e^r\qquad \text{for\ all}\ r>0.
\]
See \cite[Section~3]{Ts} or \cite[Section~2 and Appendix~A]{Suyari} for more details on $q$-logarithmic, exponential functions.
The radial derivative $\nabla_r u(x)$ 
and
the tangential derivative 
$\nabla_{\mathbb{S}^{n-1}}u(x)$ are defined by
\begin{equation*}
\nabla_ru(x):=\left(\frac{x}{|x|}\cdot \nabla u(x)\right)\frac{x}{|x|},\qquad 
\nabla_{\mathbb{S}^{n-1}}u(x):=\nabla u(x)-\nabla_r u(x).
\end{equation*}
We define the differential operator $L_{p,\theta}$ as follows:
\begin{equation*}
L_{p,\theta}u
:=
\nabla_{\mathbb{S}^{n-1}}u(x)
\slash
\left[
(q-1)
\log_{q}\frac{R}{|x|}
\right]
+
\nabla_r u(x).
\end{equation*}
The function of the Aubin--Talenti type~\eqref{eq:1.9a} and the scaling~\eqref{eq:2.1a} can be written 
by using the $q$-logarithmic and exponential functions, that is,
\[
\begin{aligned}
U_R(x)
&
=
\left[a+b
\left\{
\frac{1}{|x|^{\frac{n-p}{p-1}}}
-
\frac{1}{R^{\frac{n-p}{p-1}}}
\right\}^{-\frac{p}{n-p}}\right]^{1-\frac{n}{p}}
\\
&
=
\left[a+bR^{\frac{p}{p-1}}
\left\{
-\frac{n-p}{p-1}
\log_{\frac{n-1}{p-1}}\frac{|x|}{R}
\right\}^{-\frac{p}{n-p}}\right]^{1-\frac{n}{p}}
\end{aligned}
\]
and
\[
\left\{
\begin{aligned}
 &
 u_{\lambda}(x)
 := \lambda^{-\frac{\theta-1}{\theta}}u(x_{\lambda}),
 \\
 &
 x_{\lambda}
 :=
 \left[
  \lambda |x|^{-\frac{n-p}{p}\frac{\theta}{\theta-1}}+(1-\lambda)R^{-\frac{n-p}{p}\frac{\theta}{\theta-1}}
 \right]^{-\frac{p}{n-p}\frac{\theta-1}{\theta}}\frac{x}{|x|}
=
 R\exp_{q}\left[\lambda \log_{q}\frac{|x|}{R}\right]\frac{x}{|x|}.
\end{aligned}
\right.
\]
These expressions by $q$-functions will be used in the proof of Theorem~\ref{Theorem:1.1}.

\section{Proof of Theorems
}
Two proofs of Theorem~\ref{Theorem:1.1} will be given.
One proof is based on the scale invariance under \eqref{eq:2.1a}, 
and for the other one,
we prove 
the equivalence between the classical inequalities~\eqref{eq:1.5bbb}--\eqref{eq:1.5ccc} 
and Theorem~\ref{Theorem:1.1}.
Both proofs stand on the scaling~\eqref{eq:2.1a}, therefore, we start by
their formal derivations.

\subsection{Formal derivation of the scaling~\eqref{eq:2.1a}}
Fix $\lambda>0$ and the radially symmetric functions $u,\phi_{\lambda} \in C_0^{\infty}(B_R)$.
Define
\[
 u_{\lambda}(x):=\lambda u(\phi_{\lambda}(x)),
\]
and consider the condition on $\phi_{\lambda}$ so that
\[
\int_{B_R}|x|^{\frac{n}{p}\theta}|\nabla u_{\lambda}(x)|^{\theta}\frac{dx}{|x|^n}
=
\int_{B_R}|x|^{\frac{n}{p}\theta}|\nabla u(x)|^{\theta}\frac{dx}{|x|^n}.
\]
Since $u$ and $\phi_{\lambda}$ are radially symmetric, we have
\[
\int_{B_R}|x|^{\frac{n}{p}\theta}|\nabla u_{\lambda}(x)|^{\theta}\frac{dx}{|x|^n}
=
\omega_{n-1}
\int_0^R
\lambda^\theta
 r^{\frac{n}{p}\theta-1} \left|\frac{d}{dr}u_{\lambda}(\phi_{\lambda}(r))\right|^{\theta} dr.
\]
Changing the variable $s=\phi_{\lambda}(r)$, we have
\[
\begin{aligned}
\int_0^R
\lambda^\theta
 r^{\frac{n}{p}\theta-1} 
 \left|\frac{d}{dr}u_{\lambda}(\phi_{\lambda}(r))\right|^{\theta} dr
&
=
\int_0^R 
\lambda^{\theta}
 r^{\frac{n}{p}\theta-1}
\left|s' u'(s)\right|^{\theta} \frac{1}{s'} ds
\\
&
=
\int_0^R 
\lambda^{\theta}
 r^{\frac{n}{p}\theta-1}s'^{\theta-1}
\left|u'(s)\right|^{\theta}ds.
\end{aligned}
\]
Therefore, if $s=\phi_{\lambda}(r)$ satisfies 
the ordinary differential equation
\begin{equation}\label{eq:A1}
\lambda^{\theta}
 r^{\frac{n}{p}\theta-1}\phi_{\lambda}'^{\theta-1}
=
\phi_{\lambda}^{\frac{n}{p}\theta-1},
\end{equation}
we obtain the desired invariance.
Here \eqref{eq:A1} is solvable by separation of variables
and its solution under the boundary condition $\phi_{\lambda}(R)=R$ is
\[
\phi_{\lambda}(r)=\left(\lambda^{-\frac{\theta}{\theta-1}}r^{-\frac{n-p}{p}\frac{\theta}{\theta-1}}+\left(1-\lambda^{-\frac{\theta}{\theta-1}}\right)R^{-\frac{n-p}{p}\frac{\theta}{\theta-1}}\right)^{-\frac{p}{n-p}\frac{\theta-1}{\theta}}.
\]
Taking $\lambda \mapsto \lambda^{-\frac{\theta-1}{\theta}}$,
we obtain the scaling~\eqref{eq:2.1a}.

\subsection{Proof of Theorem~\ref{Theorem:1.1} : Scale invariance approach}
%
%
%


We first prove the following lemma:
\begin{Lemma}\label{Lemma:2.1a}
Let $\phi\in C^1(\mathbb{R})$ and $x=\phi(|y|)\frac{y}{|y|}$.
Then, the following must  hold
\begin{equation}\label{eq:2.11a}
\frac{\partial x_i}{\partial y_j}
=
\frac{|x|}{|y|}
\left[
\delta_{ij}+
\left(
\frac{|y|}{|x|}\phi'(|y|)-1
\right)
\frac{y_i y_j}{|y|^2}
\right],
\end{equation}
\begin{equation}\label{eq:2.12a}
\det\left(\frac{\partial x}{\partial y}\right)
=
\frac{|y|}{|x|}\phi'(|y|)\frac{|x|^n}{|y|^n},
\end{equation}
and
\begin{equation}\label{eq:2.13a}
\left|
\nabla_y u(x)
\right|^2
=
\left|
\frac{\partial x}{\partial y}\nabla_x u(x)
\right|^2
=
\frac{|x|^2}{|y|^2}
\left(
\left|\nabla_{\mathbb{S}^{n-1}}u(x)\right|^2
+
\left(
\frac{|y|}{|x|}\phi'(|y|)
\right)^2
\left|
\nabla_r u(x)
\right|^2
\right).
\end{equation}
\end{Lemma}

\begin{prf}{Lemma~\ref{Lemma:2.1a}}
The equality~\eqref{eq:2.11a} can be obtained by direct computations. 
Let $a:=\frac{|y|}{|x|}\phi'(|y|)-1$.
By 
\eqref{eq:2.11a},
the Jacobian is given by
\[
\det \left(\frac{\partial x}{\partial y}\right)
 =
 \frac{|x|^n}{|y|^n}
 \det
 \left(
    \begin{array}{cccc}
 a\frac{y_1^2}{|y|^2}+1 & a\frac{y_1 y_2}{|y|^2} & \cdots & a\frac{y_1y_n}{|y|^2}\\
 a\frac{y_2y_1}{|y|^2} & a\frac{y_2^2}{|y|^2}+1 & \cdots & \vdots\\
 \vdots & \vdots & \ddots & \vdots\\
 a\frac{y_ny_1}{|y|^2} & \cdots &\cdots & a\frac{y_n^2}{|y|^2}+1
    \end{array}
  \right).
\]
By the standard argument, one can see that $1$ and $(a+1)$ are eigenvalues of the matrix,
and
their multiplicities are $n-1$ and $1$, respectively. 
Therefore,
\begin{equation*}
 \det
 \left(
    \begin{array}{cccc}
 a\frac{y_1^2}{|y|^2}+1 & a\frac{y_1 y_2}{|y|^2} & \cdots & a\frac{y_1y_n}{|y|^2}\\
 a\frac{y_2y_1}{|y|^2} & a\frac{y_2^2}{|y|^2}+1 & \cdots & \vdots\\
 \vdots & \vdots & \ddots & \vdots\\
 a\frac{y_ny_1}{|y|^2} & \cdots &\cdots & a\frac{y_n^2}{|y|^2}+1
    \end{array}
  \right)
 =(a+1)
=\frac{|y|}{|x|}\phi'(|y|).
 \end{equation*}
This proves~\eqref{eq:2.12a}.
Now we need to prove \eqref{eq:2.13a}.
Since
\begin{equation*}
\begin{aligned}
\frac{\partial x}{\partial y}\nabla u(x)
=
\left(
\sum_{j=1}\frac{\partial x_i}{\partial y_j}\partial_{x_j}u
\right)_{i=1,\ldots, n}
=
\frac{|x|}{|y|}
\left(
\nabla u+
a \frac{x}{|x|}
\nabla_r u
\right),
\end{aligned}
\end{equation*}
we have
\begin{equation*}
\begin{aligned}
\left|
\frac{\partial x}{\partial y}\nabla u(x)
\right|^2
&
=
\frac{|x|^2}{|y|^2}
\Bigl(
|\nabla u|^2
+
2a|\nabla_r u|^2
+a^2|\nabla_r u|^2
\Bigr)
\\
&
=
\frac{|x|^2}{|y|^2}
\Bigl(
|\nabla_{\mathbb{S}^{n-1}}u|^2
+(a+1)^2|\nabla_r u|^2
\Bigr).
\end{aligned}
\end{equation*}
This yields \eqref{eq:2.13a} and completes the proof of Lemma~\ref{Lemma:2.1a}.
\end{prf}


\bigskip

\begin{prf}{Theorem~\ref{Theorem:1.1}(The inequality)}
We start the proof from the classical Caffarelli--Kohn--Nirenberg inequality. 
Let $1< p<n,\ \theta> 1,\ 0\le \frac{1}{\theta}-\frac{1}{\sigma}\le \frac{1}{n}$ and 
$u\in C_0^{\infty}(B_R)$. Recall that
the rescaled function $u_{\lambda}$ defined in~\eqref{eq:1.4} belongs to $C_0^{\infty}(B_R)$.
Then we have
\begin{equation}\label{eq:1.5}
S
\left(
 \int_{B_R}|x|^{\frac{n\sigma}{p^*}}|u_{\lambda}(x)|^{\sigma}\frac{dx}{|x|^n}
\right)^{\frac{1}{\sigma}}
\le
\left(
 \int_{B_R}|x|^{\frac{n\theta}{p}}|\nabla u_{\lambda}(x)|^{\theta}\frac{dx}{|x|^n}
\right)^{\frac{1}{\theta}}.
\end{equation}
We calculate both the sides of \eqref{eq:1.5}
 by considering the change of variables 
\[
\begin{aligned}
z=x_{\lambda}
&
=
 \left[
  \lambda |x|^{-\frac{n-p}{p}\frac{\theta}{\theta-1}}+(1-\lambda)R^{-\frac{n-p}{p}\frac{\theta}{\theta-1}}
 \right]^{-\frac{p}{n-p}\frac{\theta-1}{\theta}}\frac{x}{|x|}
\\
&
=
R\exp_{q}\left[\lambda \log_{q}\frac{|x|}{R}\right]\frac{x}{|x|}.
\end{aligned}
\]
Let $\phi(r)=R\exp_{q}\left[\lambda \log_{q}\frac{r}{R}\right]$.
Then,
\[
\phi'(r)
=
R\left(\exp_{q}\left[\lambda \log_{q}\frac{r}{R}\right]\right)^{q}\lambda \left(\frac{r}{R}\right)^{-q}\frac{1}{R}
=
\lambda \left(\frac{|z|}{|x|}\right)^{q}.
\]
Hence, we have
\begin{equation}\label{eq:3.21c}
\frac{|x|}{|z|}\phi'(|x|)
=
\lambda\left(\frac{|x|}{|z|}\right)^{1-q}
=
1-
\left(
\frac{R}{|z|}
\right)^{1-q}
+\lambda
\left(
\frac{R}{|z|}
\right)^{1-q}
=
(q-1)\log_{q}\frac{R}{|z|}
+\lambda
\left(
\frac{R}{|z|}
\right)^{1-q}
\end{equation}
and
\begin{equation}\label{eq:3.22c}
|x|^{\frac{n-p}{p}}\lambda^{-\frac{\theta-1}{\theta}}
=
|z|^{\frac{n-p}{p}}
\left[(q-1)\log_{q}\frac{R}{|z|}+\lambda\left(\frac{R}{|z|}\right)^{1-q}\right]^{-\frac{\theta-1}{\theta}}.
\end{equation}
From \eqref{eq:3.21c}, \eqref{eq:3.22c}, and Lemma~\ref{Lemma:2.1a}, we have
\begin{equation}\label{eq:2.4}
\begin{aligned}
 \int_{B_R}|x|^{\frac{n\sigma}{p^*}}|u_{\lambda}(x)|^{\sigma}\frac{dx}{|x|^n}
& =
 \int_{B_R}
  |x|^{\frac{n-p}{p}\sigma-n}
   {\lambda}^{-\frac{\theta-1}{\theta}\sigma}
  |u(z)|^{\sigma}
  \left|\det\left(\frac{\partial x}{\partial z}\right)\right|
 dz
\\
& =
  \int_{B_R}
  |z|^{\frac{n}{p^*}\sigma-n}
  \frac{
   |u(z)|^{\sigma}
  }
  {
   \left[
   \displaystyle 
\lambda \left(\frac{|z|}{R}\right)^{q-1}
+
\left(q-1\right)\log_{q}\frac{R}{|z|}
\right]^{1+\frac{\theta-1}{\theta}\sigma}
  }
 dz.
\end{aligned}
\end{equation}
The monotone convergence theorem shows us that the right hand side of \eqref{eq:2.4} converges to the right hand side of Theorem~\ref{Theorem:1.1}
as $\lambda\to 0$.
Similarly, Lemma~\ref{Lemma:2.1a} along with \eqref{eq:3.21c} and \eqref{eq:3.22c} implies that
\begin{equation*}
\begin{aligned}
 \int_{B_R}
&
|x|^{\frac{n\theta}{p}}|\nabla u_{\lambda}(x)|^{\theta}\frac{dx}{|x|^n}
\\
&=
 \int_{B_R}
 |x|^{\frac{n\theta}{p}}
 \lambda^{-\frac{\theta-1}{\theta}\theta}
 \left|\left(\frac{\partial z}{\partial x}\nabla u\right)(z)\right|^{\theta}
 \left|\det\left(\frac{\partial x}{\partial z}\right)\right|
 \frac{dz}{|x|^n}
 \\
&=
\int_{B_R}
|z|^{\frac{n}{p}\theta}
\frac{
\left[
   |\nabla_{\mathbb{S}^{n-1}} u|^2
   +
   \left(
\lambda \left(\frac{|z|}{R}\right)^{q-1}
+
\left(q-1\right)\log_{q}\frac{R}{|z|}
   \right)^2 \left|\nabla_r u\right|^2
  \right]^{\frac{\theta}{2}}
}
{
   \left[
\lambda \left(\frac{|z|}{R}\right)^{q-1}
+
\left(q-1\right)\log_{q}\frac{R}{|z|}
\right]^{\theta}
}
\frac{dz}{|z|^n}.
 \end{aligned}
\end{equation*}
The dominated convergence theorem as $\lambda\to 0$ yields the desired term in \eqref{eq:3.17a}.

The homogeneous inequality \eqref{eq:3.17ccc} for $\theta=\sigma$ can be obtained by starting the argument
from
\eqref{eq:1.5ccc} instead of from \eqref{eq:1.5bbb}. 
Indeed, the equality
\[
 \int_{B_R}
|x|^{\frac{n\theta}{p}}|\nabla_r u_{\lambda}(x)|^{\theta}\frac{dx}{|x|^n}
=
 \int_{B_R}
|x|^{\frac{n\theta}{p}}|\nabla_r u(x)|^{\theta}\frac{dx}{|x|^n}
\]
and 
\eqref{eq:2.4}
yield 
the homogeneous inequality \eqref{eq:3.17ccc} by taking $\lambda \to 0$.
This completes the proof of the inequalities \eqref{eq:3.17a} and \eqref{eq:3.17ccc}.
%
\end{prf}

It now remains to be proven that \eqref{eq:3.17a} and \eqref{eq:3.17ccc} are scale invariant under the scaling \eqref{eq:2.1a}. 
This can be done by direct calculation; however, we prove it in the next subsection by focusing the equivalence 
between classical inequalities~\eqref{eq:1.5bbb}--\eqref{eq:1.5ccc} 
and Theorem~\ref{Theorem:1.1}.
We prove this at the end of Section~\ref{Section:2.2}.




\subsection{Proof of Theorem~\ref{Theorem:1.1} : Transformation}
\label{Section:2.2}
The following equivalence 
gives Theorem~\ref{Theorem:1.1}. 
\begin{Proposition}\label{Proposition:3.1}
Let $p<n,\ \theta> 1,\ 0\le \frac{1}{\theta}-\frac{1}{\sigma}\le \frac{1}{n}$,
$\frac{1}{p^{*}}=\frac{1}{p}-\frac{1}{n}$
and
$q=1+\frac{n-p}{p}\frac{\theta}{\theta-1}$.
Then the following equalities 
\begin{equation}\label{eq:3.18a}
\begin{aligned}
 \int_{\mathbb{R}^n}
 |y|^{\frac{n\sigma}{p^{*}}}|v(y)|^{{\sigma}}
 \frac{dy}{|y|^{{n}}}
&
=
 \int_{{B_{R}}}\frac{|x|^{\frac{n\sigma}{p^{*}}}|u(x)|^{{\sigma}}
 }
 {
 \left[
(q-1)
 \log_{{q}}\frac{R}{|x|}\right]^{{1+\frac{\theta-1}{\theta}\sigma}}
 }
 \frac{dx}{|x|^{{n}}},
\\
 \int_{\mathbb{R}^n}
 |y|^{\frac{n\theta}{p}}|\nabla v(y)|^{{\theta}}
 \frac{dy}{|y|^{{n}}}
&
=
 \int_{{B_{R}}}
 |x|^{\frac{n\theta}{p}}
 |L_{p,\theta} u(x)|^{\theta}
 \frac{dx}{|x|^{{n}}},
\\
 \int_{\mathbb{R}^n}
 |y|^{\frac{n\theta}{p}}|\nabla_r v(y)|^{{\theta}}
 \frac{dy}{|y|^{{n}}}
&
=
 \int_{{B_{R}}}
 |x|^{\frac{n\theta}{p}}
 |\nabla_r u(x)|^{\theta}
 \frac{dx}{|x|^{{n}}}
\end{aligned}
\end{equation}
hold under the transformation
\begin{equation*}
u(x)=v(y),\qquad y=R\left(
-(q-1)\log_{q}\frac{|x|}{R}\right)^{-\frac{1}{q-1}}\frac{x}{|x|}.
\end{equation*}
\end{Proposition}


\begin{prf}{Proposition~\ref{Proposition:3.1}}
For parameters $\alpha>0$ and $q>1$, 
we consider the transformation
\begin{equation*}
 u(x)=v(y)
 ,\qquad y=C\left(-\log_{q}\frac{|x|}{R}\right)^{-\alpha}\frac{x}{|x|}
\end{equation*}
or equivalently
\[ 
x=R\exp_{q}\left[-\left(\frac{|y|}{C}\right)^{-\frac{1}{\alpha}}\right]\frac{y}{|y|}
\]
and choose suitable $\alpha,q$, and $C$ later.
It follows from
\[
-\log_{q}\frac{|x|}{R}
=
-\frac{\left(\frac{|x|}{R}\right)^{1-q}-1}{1-q}
=
\left(\frac{|x|}{R}\right)^{1-q}
\frac{\left(\frac{R}{|x|}\right)^{1-q}-1}{1-q}
=
\left(\frac{|x|}{R}\right)^{1-q}
\left(\log_{q}\frac{R}{|x|}\right)
\]
that
\begin{equation}
\label{eq:3.3a}
y
=C\left(-\log_{q}\frac{|x|}{R}\right)^{-\alpha}\frac{x}{|x|}
=
C\left(\frac{|x|}{R}\right)^{\alpha(q-1)}\left(\log_{q}\frac{R}{|x|}\right)^{-\alpha}\frac{x}{|x|}.
\end{equation}
Let $\phi(r)=R\exp_{q}[-\left(\frac{r}{C}\right)^{-\frac{1}{\alpha}}]$. 
Direct computations with \eqref{eq:3.3a} give us the following:
\begin{equation}\label{eq:2.21a}
\begin{aligned}
\frac{|y|}{|x|}\phi'(|y|)
&
=
\frac{|y|}{|x|}
R
\left(
\exp_{q}\left[-\left(\frac{|y|}{C}\right)^{-\frac{1}{\alpha}}\right]
\right)^{q}
\frac{1}{\alpha C}
\left(\frac{|y|}{C}\right)^{-\frac{1}{\alpha}-1}
\\
&
=
\frac{1}{\alpha }
\left(
\frac{|x|}{R}
\right)^{q-1}
\left(\frac{|y|}{C}\right)^{-\frac{1}{\alpha}}
=
\frac{1}{\alpha}\log_{q}\frac{R}{|x|}.
\end{aligned}
\end{equation}
Now, we are ready to prove 
\eqref{eq:3.18a}.
By applying 
Lemma~\ref{Lemma:2.1a} with the help of \eqref{eq:3.3a} and \eqref{eq:2.21a}, we have
\begin{equation}\label{eq:3.40ccc}
\begin{aligned}
 \int_{\mathbb{R}^n}
&
 |y|^{\frac{n\theta}{p}}|\nabla v(y)|^{{\theta}}
 \frac{dy}{|y|^{{n}}}
 \\
&
=
 \int_{B_R}
 |y|^{\frac{n\theta}{p}-n}
 \left|
 \frac{\partial x}{\partial y}\nabla u(x)
 \right|^{{\theta}}
 \left|
\det\left(\frac{\partial y}{\partial x}\right)
 \right|
dx
\\
&
=
 \int_{B_R}
 |y|^{\frac{n\theta}{p}-n}
 \frac{|x|^{\theta}}{|y|^{\theta}}
\left(
|\nabla_{\mathbb{S}^{n-1}}u|^2
+
  \left(
   \frac{1}{\alpha}\log_{q}\frac{R}{|x|}
  \right)^2|\nabla_r u|^2
\right)^{\frac{\theta}{2}}
 \frac{|y|^n}{|x|^n}\frac{\alpha}{\log_{q}\frac{R}{|x|}}dx
\\
&
=
\alpha^{-\alpha\theta\frac{n-p}{p}} R^{-\alpha(q-1)\frac{n-p}{p}\theta}C^{\frac{n-p}{p}\theta}
\\
&
\qquad 
\times
 \int_{{B_{R}}}
 \frac{
 |x|^{\alpha(q-1)\frac{n-p}{p}\theta+\theta}
 \left(
|\nabla_{\mathbb{S}^{n-1}}u(x)|^2
+
\left[\frac{1}{\alpha}\log_{{q}}\frac{R}{|x|}\right]^{2}
|\nabla_r u(x)|^2
 \right)^{{\frac{\theta}{2}}}
 }
 {
 \left(\frac{1}{\alpha}\log_{{q}}\frac{R}{|x|}\right)^{\alpha \theta \frac{n-p}{p}+1}
 }
 \frac{dx}{|x|^{{n}}}.
\end{aligned}
\end{equation}
Here we take 
\begin{equation}\label{eq:2.30ccc}
\alpha=\frac{p}{n-p}\frac{\theta-1}{\theta},\qquad q=1+\frac{n-p}{p}\frac{\theta}{\theta-1},
\qquad
C=R\alpha^{\alpha}
\end{equation}
so that 
\[
\begin{aligned}
&
\alpha(q-1)\theta\frac{n-p}{p}+\theta=\frac{n\theta}{p},
\quad
&&
\alpha \theta \frac{n-p}{p}+1=\theta,
\\
&
\alpha^{-\alpha\theta\frac{n-p}{p}} R^{-\alpha(q-1)\frac{n-p}{p}\theta}C^{\frac{n-p}{p}\theta}=1,
\quad 
&&
\frac{1}{\alpha}=q-1.
\end{aligned}
\]
Therefore, the last term in \eqref{eq:3.40ccc} coincides with 
$
\displaystyle
 \int_{{B_{R}}}
 |x|^{\frac{n\theta}{p}}
|L_{p,\theta} u(x)|^{\theta}
 \frac{dx}{|x|^{{n}}}.
$
Applying a similar argument, it follows that
\[
 \int_{\mathbb{R}^n}
 |y|^{\frac{n\theta}{p}}|\nabla_r v(y)|^{{\theta}}
 \frac{dy}{|y|^{{n}}}
=
 \int_{{B_{R}}}
 |x|^{\frac{n\theta}{p}}
 |\nabla_r u(x)|^{\theta}
 \frac{dx}{|x|^{{n}}}.
\]
Finally, \eqref{eq:2.12a}, \eqref{eq:2.21a}, and \eqref{eq:3.3a} 
imply that
\begin{equation*}
\begin{aligned}
 \int_{\mathbb{R}^n}
 |y|^{\frac{n\sigma}{p^{*}}}|v(y)|^{{\sigma}}
 \frac{dy}{|y|^{{n}}}
&
=
 \int_{B_R}
 |y|^{\frac{n\sigma}{p^{*}}-n}|u(x)|^{{\sigma}}
  \left|
\det\left(\frac{\partial y}{\partial x}\right)
 \right|
 dx
\\
&
=
 \int_{B_R}
 |y|^{\frac{n\sigma}{p^{*}}-n}|u(x)|^{{\sigma}}
\frac{|y|^n}{|x|^n}
\frac{\alpha}{\log_{q}\frac{R}{|x|}}
 dx
\\
&
=
\alpha^{-\frac{\theta-1}{\theta}\sigma} R^{-\frac{n-p}{p}\sigma}C^{\frac{n-p}{p}\sigma}
 \int_{{B_{R}}}
 \frac{|x|^{\frac{n\sigma}{p^{*}}}|u(x)|^{{\sigma}}
 }
 {
 \left[\frac{1}{\alpha}\log_{{q}}\frac{R}{|x|}\right]^{{1+\frac{\theta-1}{\theta}\sigma}}
 }
 \frac{dx}{|x|^{{n}}}.
\end{aligned}
\end{equation*}
Applying 
\eqref{eq:2.30ccc},
we have
$\alpha^{-\frac{\theta-1}{\theta}\sigma} R^{-\frac{n-p}{p}\sigma}C^{\frac{n-p}{p}\sigma}=1$ and $\frac{1}{\alpha}=q-1$.
This proves Proposition~\ref{Proposition:3.1}.
\end{prf}

As a consequence of Proposition~\ref{Proposition:3.1},
we obtain the inequalities~\eqref{eq:3.17a} and \eqref{eq:3.17ccc} from the classical inequalities~\eqref{eq:1.5bbb}
and \eqref{eq:1.5ccc}, respectively.
Moreover, the function of the Aubin--Talenti type~\eqref{eq:1.9a} in Theorems~\ref{Theorem:1.2}--\ref{Theorem:1.3} is obtained
from 
the original Aubin--Talenti function~\eqref{AubinTalenti} in $\mathbb{R}^n$ via the transformation.
%
%

Finally,
we give a simple proof of the scale invariance of the inequalities~\eqref{eq:3.17a}-\eqref{eq:3.17ccc} by focusing the scaling structure
between the dilation~\eqref{eq:1.7a} and the scaling~\eqref{eq:2.1a}.
Let $u\in C_0^{\infty}(B_R)$, $v\in C_0^{\infty}(\mathbb{R}^n)$, $\alpha>0$,
$q>1$, and $\lambda>0$.
We define the following transformations:
\begin{equation*}
\begin{aligned}
&
D_{\lambda}:C_0^{\infty}(\mathbb{R}^n)\to C_0^{\infty}(\mathbb{R}^n);
\quad 
D_{\lambda}v(y):=v(\lambda y),
\\
&
S_{\lambda}:C_0^{\infty}(B_R)\to C_0^{\infty}(B_R);
\quad 
S_{\lambda}u(x):=u(x_{\lambda}),
\quad
x_{\lambda}:=R\exp_{q}\left[\lambda \log_{q}\frac{|x|}{R}\right]
\frac{x}{|x|}
\\
&
T_{\alpha}:C_0^{\infty}(\mathbb{R}^n) \to C_0^{\infty}(B_R);
\quad
T_{\alpha,q}v(x):=v\left(y\right),
\quad y:=
R\left(-(q-1)\log_{q}\frac{|x|}{R}\right)^{-\alpha} \frac{x}{|x|}
\end{aligned}
\end{equation*}
and its inverse transformation $T_{\alpha,q}^{-1}$.

%

\begin{Theorem}\label{Theorem:3.2a}
For all $\lambda,\ \alpha>0$, $q>1$, $u\in C_{0}^{\infty}(B_R)$, and $v\in C_0^{\infty}(\mathbb{R}^n)$,
there hold
\[
\begin{aligned}
 S_{\lambda}(T_{\alpha,q}v)(x)
 &
 =
 T_{\alpha}(D_{\lambda^{-\alpha}}v)(x)\quad \text{for\ all\ }x\in B_R,
 \\
 T_{\alpha}^{-1}(S_{\lambda}u)(y)
 &
 =
 D_{\lambda^{-\alpha}}(T^{-1}_{\alpha,q}u)(y)\quad \text{for\ all\ }y\in \mathbb{R}^n.
\end{aligned}
\]
\end{Theorem}
\begin{prf}{Theorem~\ref{Theorem:3.2a}}
Fix $v\in C_0^{\infty}(\mathbb{R}^n)$. 
Then, 
\[
\begin{aligned}
 S_{\lambda}(T_{\alpha,q}v)(x)
&
=
v\left(
\left(
 -\lambda (q-1) \log_{q}\frac{|x|}{R}
\right)^{-\alpha}
\right)
\\
&=
(D_{\lambda^{-\alpha}}v)\left(
\left(
 -(q-1)\log_{q}\frac{|x|}{R}
\right)^{-\alpha}
\right)
=
T_{\alpha,q}(D_{\lambda^{-\alpha}}v)(x)
\end{aligned}
\]
for all $x\in B_R$.
The other equality is obtained by taking $v=T_{\alpha,q}^{-1}u$.
\end{prf}

\bigskip

Applying Theorem~\ref{Theorem:3.2a} and Proposition~\ref{Proposition:3.1} with the help of the dilation invariance of the classical inequalities~\eqref{eq:1.5bbb} and \eqref{eq:1.5ccc},
we prove the scale invariance properties of \eqref{eq:3.17a} and \eqref{eq:3.17ccc} under the scaling~\eqref{eq:2.1a}.

\bigskip

\begin{prf}{Theorem~\ref{Theorem:1.1} (Scale invariance)}
Let us define some notations for norms in Theorem~\ref{Theorem:1.1}, which simplify the proof of the scale invariance. 
\begin{equation}\label{eq:3.18a}
\begin{aligned}
&
\|v\|_{X(\mathbb{R}^n)}
:=
\left(
 \int_{\mathbb{R}^n}
 |y|^{\frac{n\sigma}{p^{*}}}|v(y)|^{{\sigma}}
 \frac{dy}{|y|^{{n}}}
\right)^{{\frac{1}{\sigma}}},
&&
\|v\|_{Y(\mathbb{R}^n)}
:=
\left(
 \int_{\mathbb{R}^n}
 |y|^{\frac{n\theta}{p}}|v(y)|^{{\theta}}
 \frac{dy}{|y|^{{n}}}
\right)^{{\frac{1}{\theta}}},
\\
&
\|u\|_{\tilde{X}(B_R)}
:=
\left(
 \int_{{B_{R}}}\frac{|x|^{\frac{n\sigma}{p^{*}}}|u(x)|^{{\sigma}}
 }
 {
 \left[(q-1)\log_{q}\frac{R}{|x|}\right]^{{1+\frac{\theta-1}{\theta}\sigma}}
 }
 \frac{dx}{|x|^{{n}}}
\right)^{{\frac{1}{\sigma}}},
&&
\|u\|_{\tilde{Y}(B_R)}
:=
\left(
 \int_{{B_{R}}}
 |x|^{\frac{n\theta}{p}}
|u(x)|^{\theta}
 \frac{dx}{|x|^{{n}}}
\right)^{{\frac{1}{\theta}}}.
\end{aligned}
\end{equation}
Assume that 
$\alpha,q$ satisfy \eqref{eq:2.30ccc}.
Then, Proposition~\ref{Proposition:3.1} can be stated as
\begin{equation}\label{eq:2.55ccc}
 \|T_{\alpha,q}^{-1}u\|_{X(\mathbb{R}^n)}=
 \|u\|_{\tilde{X}(B_R)}
 \quad \text{and} \quad 
 \|\nabla (T_{\alpha,q}^{-1}u)\|_{Y(\mathbb{R}^n)}
 =
 \|L_{p,\theta} u\|_{\tilde{Y}(B_R)}.
\end{equation}
Furthermore, the scaling~\eqref{eq:2.1a} coincides with 
$u_{\lambda}(x)=\lambda^{-\frac{\theta-1}{\theta}}(S_{\lambda}u)(x)$.

Using these notation, we prove the invariance of $\|u\|_{\tilde{X}(B_R)}$
and
$\|L_{p,\theta} u\|_{\tilde{Y}(B_R)}$ under 
the scaling $\lambda^{-\frac{\theta-1}{\theta}}S_{\lambda}$.
It follows from \eqref{eq:2.55ccc}
that
\[
\begin{aligned}
&
 \|\lambda^{-\frac{\theta-1}{\theta}}S_{\lambda}u\|_{\tilde X(B_R)}
 =
 \|\lambda^{-\frac{\theta-1}{\theta}}T_{\alpha,q}^{-1}(S_{\lambda}u)\|_{X(\mathbb{R}^n)},
\\
& 
 \|\lambda^{-\frac{\theta-1}{\theta}}L_{p,\theta}(S_{\lambda} u)\|_{\tilde Y(B_R)}
=
 \left\|\lambda^{-\frac{\theta-1}{\theta}}\nabla \bigl(T_{\alpha,q}^{-1}(S_{\lambda}u)\bigr)\right\|_{Y(\mathbb{R}^n)}.
\end{aligned}
\]
Applying Theorem~\ref{Theorem:3.2a} and the dilation invariance in $\mathbb{R}^n$ with
$\lambda^{-\frac{\theta-1}{\theta}}=(\lambda^{-\alpha})^{\frac{n-p}{p}}$,
we obtain
\[
\begin{aligned}
&
\|\lambda^{-\frac{\theta-1}{\theta}}T_{\alpha,q}^{-1}(S_{\lambda}u)\|_{X(\mathbb{R}^n)}
=
\|\lambda^{-\frac{\theta-1}{\theta}}D_{\lambda^{-\alpha}}(T_{\alpha,q}^{-1}u)\|_{X(\mathbb{R}^n)}
=
\|
T_{\alpha,q}^{-1}u\|_{X(\mathbb{R}^n)},
\\
 &
 \left\|\lambda^{-\frac{\theta-1}{\theta}}\nabla \bigl(T_{\alpha,q}^{-1}(S_{\lambda}u)\bigr)\right\|_{Y(\mathbb{R}^n)}
 =
 \left\|\lambda^{-\frac{\theta-1}{\theta}}\nabla \bigl(D_{\lambda^{-\alpha}}(T_{\alpha,q}^{-1}u\bigr)\right\|_{Y(\mathbb{R}^n)}
=
 \left\|
 \nabla \bigl(
 T_{\alpha,q}^{-1}u\bigr)\right\|_{Y(\mathbb{R}^n)}.
\end{aligned}
\]
Then, the desired invariance holds by 
applying 
\eqref{eq:2.55ccc}
 to the last equalities.
\end{prf}

%

%




\bigskip

\noindent 
{\bf Acknowledgments.} 
This work was partially funded by JSPS KAKENHI \# 18K13441.

%

\end{document}